\documentclass[a4paper]{amsart}

\usepackage{latexsym}
\usepackage{amssymb,latexsym}
\usepackage{amsmath}
\usepackage{amsthm}
\usepackage{mathrsfs}
\usepackage{amscd}
\usepackage{amsfonts}
\usepackage{enumerate}

\newcommand{\Hess}{{{\rm Hess}\,}}

\renewcommand{\thetheoremName}

\newtheorem{theorem}{Theorem}[section]

\newtheorem{proposition}[theorem]{Proposition}
\newtheorem{corollary}[theorem]{Corollary}

\theoremstyle{definition}
\newtheorem{definition}[theorem]{Definition}

\newtheorem{remark}[theorem]{Remark}

\numberwithin{equation}{section}

\newcommand{\dist}{\operatorname{dist}}
\newcommand{\Vol}{\operatorname{Vol}}

\newcommand{\erre}{\mathbb{R}}
\newcommand{\Hp}{\mathbb{H}}

\newcommand{\Sup}{\operatorname{Sup}}

\def\kan{I\!\!K^n(b)}

\def\kam{I\!\!K^m(b)}
\def\Han{\mathbb{H}^n(b)}
\begin{document}

\title[Volume growth, number of ends]{Volume growth, number of ends and the topology of a complete submanifold}

\author[V. Gimeno]{Vicent Gimeno*}
\address{Departament de Matem\`{a}tiques- Institut of New Imaging Technologies,
Universitat Jaume I, Castellon, Spain.}
\email{gimenov@uji.es}
\author[V. Palmer]{Vicente Palmer*}
\address{Departament de Matem\`{a}tiques- Institut of New Imaging Technologies,
Universitat Jaume I, Castellon, Spain.}
\email{palmer@mat.uji.es}
\thanks{}
\thanks{* Work partially supported by the Caixa Castell\'{o} Foundation, and DGI grant MTM2010-21206-C02-02.}
\maketitle
\begin{abstract}
Given a complete isometric immersion $\varphi: P^m \longrightarrow N^n$ in an ambient Riemannian manifold $N^n$ with a pole and with radial sectional curvatures bounded from above by the corresponding radial sectional curvatures of a radially symmetric space $M^n_w$, we determine a set of conditions on the extrinsic curvatures of $P$ that guarantees that the immersion is proper and that $P$ has finite topology in the line of the results in \cite{Pac} and \cite{Pac2}. When the ambient manifold is a radially symmetric space, it is shown an inequality between the (extrinsic) volume growth of a complete and minimal submanifold and its number of ends which generalizes the classical inequality stated in \cite{anderson} for complete and minimal submanifolds in $\erre^n$. We obtain as a corollary the corresponding inequality between the (extrinsic) volume growth and the number of ends of a complete and minimal submanifold in the Hyperbolic space together  with Bernstein type results for such submanifolds in Euclidean and Hyperbolic spaces, in the vein of the work \cite{Kasue1}.
\end{abstract}
\section{Introduction}\label{Intr} A natural question in Riemannian geometry is to explore the influence of the  curvature conduct of a complete Riemannan  manifold on its geometric and topological properties. Classical results concernig this are the gap theorems showed by Greene and Wu in \cite{GreW2}, (see too \cite{GreW3}), and, when it is considered a minimal submanifold (properly) immersed in the Euclidean space $\erre^n$, the Berstein-type theorems showed by Anderson in \cite{anderson} and by Schoen in \cite{schoen}. Greene and Wu's results states, roughly speaking, that a Riemannian manifold with a pole and with faster than quadratic decay of its sectional curvatures is isometric to the Euclidean space. On the other hand, Anderson proved, as a corollary of a generalization of the Chern-Osserman theorem on complete and minimal submanifolds of $\erre^n$ with finite total (extrinsic) curvature, that any of such submanifolds having one end is an affine $n$-plane. More examples concerning submanifolds immersed in an ambient Riemannian manifold and the analysis of its (intrinsic and  extrinsic) curvature behavior are the gap results, (of Bernstein-type), given by Kasue and Sugahara in \cite{Kasue1} (see Theorems A and B), where an accurate (extrinsic) curvature decay forces to minimal, (or not) submanifolds with one end of the Euclidean and Hyperbolic spaces to be totally geodesic, and the gap results for minimal submanifolds in the Euclidean space with controlled scalar curvature given by Kasue in \cite{Kasue2}. 

The estimation of the number of ends of these submanifolds plays a fundamental r\^ole in all the Bernstein-type results above mentioned. In this way, it is proved in \cite{anderson} (see Theorems 4.1 and 5.1 in that paper) that given a complete and minimal submanifold $\varphi: P^m \longrightarrow \erre^n$, ($m >2$) having finite total curvature $\int_P \Vert B^P\Vert^m d\sigma <\infty$, its (extrinsic) volume growth, defined as the quotient  $\frac{\Vol(\varphi(P) \cap B^{0,n}_t)}{\omega_n t^n}$ is bounded from above by the number of ends of $P$, $\mathcal{E}(P)$, namely
\begin{equation}\label{endsvolum}
 \lim_{t \to \infty} \frac{\Vol(\varphi(P) \cap B^{0,n}_t)}{\omega_n t^n} \leq \mathcal{E}(P)
\end{equation} 
\noindent where $B^{b,n}_t$ denotes the metric $t-$ ball in the real space form of constant curvature $b$, $\kan$,  and  $\Vert B^P\Vert $ denotes the Hilbert-Schmidt norm of the second fundamental form of $P$ in $\erre^n$. If moreover $\mathcal{E}(P)=1$,  it is concluded (using inequality (\ref{endsvolum})) the Bernstein-type result above alluded, namely,  that  $P^m$ is an affine plane, i.e. totally geodesic in $\erre^n$, (see Theorem 5.2 in \cite{anderson}).

In the paper \cite{Qing} it was proved that inequality (\ref{endsvolum}) is in fact an equality when the minimal submanifold in $\erre^n$ exhibits an accurate decay of its extrinsic curvature $\Vert B^P\Vert$ and in the paper \cite{Kasue1} it was proved that, if the submanifold $P$ has only one end and the decay of its extrinsic curvature $\Vert B^P\Vert $ is faster than linear, (when the ambient space is $\erre^n$) or than exponential, (when the ambient space is $\Han$), then it is is totally geodesic.

Within this study of the behavior at infinity of complete and minimal submanifolds with finite total curvature immersed in the Euclidean space, it was proved also in \cite{anderson} and in \cite{O} that the immersion of a complete and minimal submanifold $P$ in $\erre^n$ or $\Han$ satisfying $\int_P \Vert B^P\Vert^m d\sigma <\infty$ is proper and that $P$ is of finite topological type.

We should mention here the results in \cite{Pac} and in \cite{Pac2}, where has been stated new conditions on the decay of the extrinsic curvature for a completely immersed submanifold $P$ in the Euclidean space (\cite{Pac}) and in a Cartan-Hadamard manifold (\cite{Pac2}) which guarantees the properness of the submanifold and the finiteness of its topology.

In view of these results, it seems natural to consider the following three issues: 
\begin{enumerate}

\item Can the properness/finiteness results in \cite{Pac} and \cite{Pac2} be extended to submanifolds immersed in spaces  which have not necessarily non-positive curvature?, 

\item Do we have an analogous to inequality (\ref{endsvolum}) between the extrinsic volume growth and the number of ends when we consider a minimal submanifold (properly) immersed in Hyperbolic space which exhibit an accurate extrinsic curvature decay?. 

\item Moreover, is it possible to deduce from this inequality a Bernstein-type result in the line of \cite{anderson} and \cite{Kasue1}?.
\end{enumerate}

 We provide in this paper a (partial) answer to these questions, besides other lower bounds for the number of ends for (non-minimal) submanifolds in the Euclidean and Hyperbolic spaces and other gap results related with these estimates. As a preliminary view of our results, we have the following theorems,  Theorem \ref{cor1} and Theorem \ref{cor3}, which follows directly from our Theorem \ref{minimal}. In Theorem \ref{cor1} we have the answer to the two last questions, namely, setting equation (\ref{endsvolum}), but in the Hyperbolic case, and a Bernstein-type result for minimal submanifolds in the Hyperbolic space, in the line studied by Kasue and Sugahara in \cite{Kasue1}, (see assertion (A-iv) of Theorem A). On the other hand, Theorem \ref{cor3} encompasses a slightly less general version of  assertion (A-i) of Theorem A in \cite{Kasue1}.
 
 \begin{theorem}\label{cor1}
Let $\varphi: P^m \longrightarrow \Han$ be a complete, proper and minimal immersion with $m>2$. Let us suppose that for sufficiently large $R_0$ and for  all points $x \in P$ such that $r(x)>R_0$, (i.e. outside a compact),
 $$\Vert B^P_x \Vert \leq \frac{\delta(r(x))}{e^{2\sqrt{-b}\, r(x)}}$$
where $r(x)=d_{\mathbb{H}^n(b)}(o,\varphi(x))$ is the (extrinsic) distance in $\mathbb{H}^n(b)$ of the points in $\varphi(P)$ to a fixed pole $o \in \mathbb{H}^n(b)$ such that $\varphi^{-1}(o) \neq \emptyset$ and $\delta(r)$ is a smooth function such that $\delta(r)\to 0$ when $r\to \infty$. Then:

\begin{enumerate}
\item The finite number of ends $\mathcal{E}(P)$ is related with the volume growth by
$$\Sup_{t>0}\frac{D_t(o)}{\Vol(B_t^{b,m})}\leq \mathcal{E}(P)$$
\noindent where $D_t(o)=\{x\in P : r(x)< t\}=\{x\in P : \varphi(x) \in B^{b,n}_{t}(o)\}$
is the extrinsic ball
of radius $t$ in $P$, (see Definition \ref{ExtBall}).

\item If $P$ has only one end, $P$ is totally geodesic in $\mathbb{H}^n(b)$
\end{enumerate}
\end{theorem}

When the ambient manifold is $\erre^n$, we have the following Bernstein-type result as in \cite{Kasue1}:

\begin{theorem}\label{cor3}
Let $\varphi: P^m \longrightarrow \erre^n$ be a complete non-compact, minimal and proper immersion with $m>2$. Let us suppose that for sufficiently large $R_0$ and for all points $x \in P$ such that $r(x)>R_0$, (i.e. outside the compact extrinsic ball $D_{R_0}(o)$ with $\varphi^{-1}(o) \neq \emptyset$),

 $$\Vert B^P_x \Vert \leq \frac{\epsilon(r(x))}{r(x)}$$
 
\noindent where $\epsilon(r)$ is a smooth function such that $\epsilon(r)\to 0$ when $r\to \infty$. Then:

\begin{enumerate}
\item The finite number of ends $\mathcal{E}(P)$ is related with the volume growth by
$$\Sup_{t>0}\frac{\Vol(D_t)}{\Vol(B_t^{0,m})}\leq \mathcal{E}(P)$$ 
\item If $P$ has only one end, $P$ is totally geodesic in $\erre^n$.
\end{enumerate}
\end{theorem}

These results, that we shall prove in Section \ref{proofcorollaries}, (together the corollaries of Section \ref{corollaries}), follows from two main theorems, stablished in Section \ref{Main-Results}. In the first (Theorem \ref{theorem1}) we show that a complete isometric immersion $\varphi: P^m \longrightarrow N^n$, ($m >2$), with controlled second fundamental form in a complete Riemannian manifold  which possess a pole and has controlled radial sectional curvatures is proper and has finite topology. In the second (Theorem \ref{theorem2}) it is proved that a complete and proper isometric immersion $\varphi: P^m \longrightarrow M^n_w$, ($m>2$), with controlled second fundamental form in a radially symmetric space $M^n_w$ with sectional curvatures bouded from below by a radial function has its volume growth bounded from above by a quantity which involve its (finite) number of ends.

The proof of both theorems follows basically the argumental lines of the proofs given in \cite{Pac} and \cite{Pac2} and some ideas in  \cite{Qing}. An important difference to these results is that, on our side, we allow to the ambient manifold to have positive sectional curvatures, bounding from above only the sectional curvatures of the planes containing radial directions. However, to show the properness of the immersion in \cite{Pac2}, the ambient manifold must have non-positive sectional curvatures, and to assure the finiteness of the topology of the immersion $P$, this ambient manifold must be, in addition, simply connected, (i.e. a Cartan-Hadamard manifold).  This difference is based in following considerations.

To obtain the finiteness of the topology in Theorem \ref{theorem1}, we show that the restricted, (to the submanifold) extrinsic distance to a fixed pole (in the ambient manifold) has no critical points outside a compact and then, we apply classical Morse theory. To show that the extrinsic distance function has no critical points  we compute its Hessian as we can find it in  \cite{Ma-Pa} and \cite{Pa3}. These results are, in its turn, based in the Jacobi-Index analysis for the Hessian of the distance function given in \cite{GreW}, in particular, its Theorem A, (see Subsection \ref{subsecLap}). This comparison theorem is different of the Hessian comparison Theorem 1.2 used in \cite{Pac2}: while in this last theorem, the space used as a model to compare is the real space form with constant sectional curvature equal to the bound on the sectional curvatures of the given Riemannian manifold, in our adaptation of  Theorem A in \cite{GreW}, (see Theorem \ref{thmGreW}), only the sectional curvatures of the planes containing radial directions from the pole are bounded by the corresponding radial sectional curvatures in a radially symmetric space used as a model. 

We also note at this point that although we use the definition of pole given by Greene and Wu in \cite{GreW}, (namely, the exponential must be a diffeomorphism at a pole), in fact, the comparison of the Hessians in Theorem A holds along radial geodesics from the poles defined as those points which have not conjugate points, as in \cite{Pac2}.

\subsection{Outline}
The outline of the paper is the following.
In Section \S.2 we present the definiton of extrinsic ball, together the basic facts about the Hessian comparison theory of restricted distance function we are going to use and an isoperimetric inequality for the extrinsic balls which plays an important r\^ole in the proof of Theorem \ref{theorem2} . Section \S.3 is devoted to the statement of the main results (Theorem \ref{theorem1}, Theorem \ref{theorem2} and Theorem \ref{minimal}). We shall present in Section \ref{corollaries} two lists of results based in  Theorems \ref{theorem1}, \ref{theorem2} and \ref{minimal}:  the first set of consequences is devoted to bound from above the volume growth of a submanifold by the number of its ends, in several contexts, obtaining moreover some Bernstein-type results. In the second set of corollaries are stated some compactification theorems for submanifolds in $\erre^n$, in $\Hp^n$ and in $\Hp^n \times \erre^l$. Sections \S.5, \S.6, \S.7 are devoted to the proof of Theorems \ref{theorem1}, \ref{theorem2}, and \ref{minimal}, respectively. Theorem \ref{cor1}, Theorem \ref{cor3} and the corollaries stated in Section \S.4 are proved in Section \S.8.

\section{Preliminaires}\label{Prelim}

\subsection{The extrinsic distance}

We assume throughout the paper that $\varphi: P^m \longrightarrow N^n$ is an isometric immersion of a complete non-compact Riemannian $m$-manifold $P^m$ into a complete Riemannian manifold $N^n$ with a pole $o\in N$, (this is the precise meaning we shall give to the word {\em submanifold} along the text) . Recall that a pole
is a point $o$ such that the exponential map
$$\exp_{o}\colon T_{o}N^{n} \to N^{n}$$ is a
diffeomorphism. For every $x \in N^{n}- \{o\}$ we
define $r(x) = r_o(x) = \dist_{N}(o, x)$, and this
distance is realized by the length of a unique
geodesic from $o$ to $x$, which is the {\it
radial geodesic from $o$}. We also denote by $r\vert_P$ or by $r$
the composition $r\circ \varphi: P\to \erre_{+} \cup
\{0\}$. This composition is called the
{\em{extrinsic distance function}} from $o$ in
$P^m$. The gradients of $r$ in $N$ and $r\vert_P$ in  $P$ are
denoted by $\nabla^N r$ and $\nabla^P r$,
respectively. Then we have
the following basic relation, by virtue of the identification, given any point $x\in P$, between the tangent vector fields $X \in T_xP$ and $\varphi_{*_{x}}(X) \in T_{\varphi(x)}N$
\begin{equation}\label{radiality}
\nabla^N r = \nabla^P r +(\nabla^N r)^\bot ,
\end{equation}
where $(\nabla^N r)^\bot(\varphi(x))=\nabla^\bot r(\varphi(x))$ is perpendicular to
$T_{x}P$ for all $x\in P$.

\begin{definition}\label{ExtBall}
Given $\varphi: P^m \longrightarrow N^n$ an isometric immersion of a complete and connected Riemannian $m$-manifold $P^m$ into a complete Riemannian manifold $N^n$ with a pole $o\in N$,  we denote the {\em{extrinsic metric balls}} of radius $t >0$ and center $o \in N$ by
$D_t(o)$. They are defined as  the subset of $P$:
$$
D_t(o)=\{x\in P : r(\varphi(x))< t\}=\{x\in P : \varphi(x) \in B^N_{t}(o)\}
$$
where $B^N_{t}(o)$ denotes the open geodesic ball
of radius $t$ centered at the pole $o$ in
$N^{n}$. Note that the set $\varphi^{-1}(o)$ can be the empty set.
\end{definition}

\begin{remark}\label{theRemk0}
When the imersion $\varphi$ is proper, the extrinsic domains $D_t(o)$
are precompact sets, with smooth boundary $\partial D_t(o)$.  The assumption on the smoothness of
$\partial D_{t}(o)$ makes no restriction. Indeed, 
the distance function $r$ is smooth in $N - \{o\}$ 
since $N$ is assumed to possess a pole $o\in N$. Hence
the composition $r\vert_P$ is smooth in $P$ and consequently the
radii $t$ that produce smooth boundaries
$\partial D_{t}(o)$ are dense in $\mathbb{R}$ by
Sard's theorem and the Regular Level Set Theorem.

\end{remark}

We now present the curvature restrictions which constitute the geometric framework of our study.

\begin{definition}
Let $o$ be a point in a Riemannian manifold $N$
and let $x \in N-\{ o \}$. The sectional
curvature $K_{N}(\sigma_{x})$ of the two-plane
$\sigma_{x} \in T_{x}N$ is then called a
\textit{$o$-radial sectional curvature} of $N$ at
$x$ if $\sigma_{x}$ contains the tangent vector
to a minimal geodesic from $o$ to $x$. We denote
these curvatures by $K_{o, N}(\sigma_{x})$.
\end{definition}


\subsection{Model spaces}\label{subsecWarp}
Throughout this paper we shall assume that the ambient manifold
$N^n$ has its $o$-radial sectional curvatures $K_{o,N}(x)$ bounded
from above by the expression $K_w(r(x))=-w''(r(x))/w(r(x))$, which are
precisely the radial sectional curvatures of the {\em
$w$-model space} $\,M^{m}_{w}\,$ we are going to define. 

\begin{definition}[See  \cite{O'N}, \cite{Gri} and  \cite{GreW}]\label{defModel}
A $w-$model $M_{w}^{m}$ is a smooth warped product with base $B^{1}
= [0,\Lambda[ \,\subset \mathbb{R}$ (where $0 < \Lambda
\leq  \infty$), fiber $F^{m-1} = \mathbb{S}^{m-1}_{1}$ (i.e. the unit
$(m-1)$-sphere with standard metric), and warping function $w\colon
[0,\Lambda[ \to \mathbb{R}_{+}\cup \{0\}$, with $w(0) = 0$,
$w'(0) = 1$, and $w(r) > 0$ for all $r >  0$. The point
$o_{w} = \pi^{-1}(0)$, where $\pi$ denotes the projection onto
$B^1$, is called the {\em{center point}} of the model space. If
$\Lambda = \infty$, then $o_{w}$ is a pole of $M_{w}^{m}$.
\end{definition}

\begin{proposition}\label{propSpaceForm}
The simply connected space forms $\mathbb{K}^{m}(b)$ of constant
curvature $b$ are $w-$models with warping functions
\begin{equation*}
w_b(r) = \begin{cases} \frac{1}{\sqrt{b}}\sin(\sqrt{b}\, r) &\text{if $b>0$}\\
\phantom{\frac{1}{\sqrt{b}}} r &\text{if $b=0$}\\
\frac{1}{\sqrt{-b}}\sinh(\sqrt{-b}\,r) &\text{if $b<0$}.
\end{cases}
\end{equation*}
Note that for $b > 0$ the function $Q_{b}(r)$ admits a smooth
extension to  $r = \pi/\sqrt{b}$.
\end{proposition}

\begin{proposition}[See Proposition 42 in Chapter 7 of \cite{O'N}. See also \cite{GreW} and \cite{Gri}]\label{propWarpMean}
Let $M_{w}^{m}$ be a $w-$model with warping function $w(r)$ and
center $o_{w}$. The distance sphere $S^w_r$ of radius $r$ and center $o_{w}$
in $M_{w}^{m}$ is the fiber $\pi^{-1}(r)$. This distance sphere has
the constant mean curvature $\eta_{w}(r)= \frac{w'(r)}{w(r)}$.

 On the other hand, the
$o_{w}$-radial sectional curvatures of $M_{w}^{m}$ at every $x \in
\pi^{-1}(r)$ (for $r > 0$) are all identical and determined
by
\begin{equation*}
K_{o_{w} , M_{w}}(\sigma_{x}) = -\frac{w''(r)}{w(r)}.
\end{equation*}

\noindent and the sectional curvatures of $M_{w}^{m}$ at every $x \in
\pi^{-1}(r)$ (for $r > 0$) of the tangent planes to the fiber $S^w_r$ are also all identical and determined by

\begin{equation*}
K(r)=K_{M_{w}}(\Pi_{S^w_r}) = \frac{1-(w'(r))^2}{w^2(r)}.
\end{equation*}

\end{proposition}

\begin{remark}
The $w-$model spaces are completely determined via $w$
by the mean curvatures of the spherical fibers $S^{w}_{r}$:
$$
\,\eta_{w}(r) = w'(r)/w(r)\,\quad,
$$
by the
volume of the fiber
$$
\,\Vol(S^{w}_{r}) \, = V_{0}\,w^{m-1}(r)\,\quad ,
$$
and by the volume of the corresponding ball, for which the fiber is
the boundary
$$
\,\Vol(B^{w}_{r}) \, = \, V_{0}\, \int_{0}^{r}
\,w^{m-1}(t)\,dt\,\quad .
$$
Here $V_{0}$ denotes the volume of the unit sphere $S^{0,m-1}_{1}$, (we denote in general as $S^{b,m-1}_r$ the sphere of radius $r$ in the real space form $\kam$) .
The latter two functions define the isoperimetric quotient function
as follows
$$
\,q_{w}(r) \, = \, \Vol(B^{w}_{r})/\Vol(S^{w}_{r}) \quad .
$$

\end{remark}

Besides the r\^ole of comparison controllers for the radial sectional curvatures of
$N^{n}$, we shall need two further purely intrinsic conditions
on the model spaces: 

\begin{definition}
A given $w-$model space $\, M^{m}_{w}\,$ is
called balanced from below and balanced from above, respectively, if
the following weighted isoperimetric conditions are satisfied:
$$
\begin{aligned}
\text{Balance from below:}\quad q_{w}(r)\,\eta_{w}(r) \, &\geq 1/m \quad \text{for all} \quad r \geq 0 \quad ;\\
\text{Balance from above:}\quad  q_{w}(r)\,\eta_{w}(r) \, &\leq
1/(m-1) \quad \text{for all} \quad r \geq 0 \quad .
\end{aligned}
$$
A model space is called {\em{totally balanced}} if it is balanced
both from below and from above. 
\end{definition}

\begin{remark}\label{propW4}
If $\,K_{w}(r) \geq -\eta_{w}^{2}(r)\,$ then $\,M_{w}^{m}\,$ is
balanced from above. If $\,K_{w}(r) \leq 0\,$ then $\,M_{w}^{m}\,$
is balanced from below, see the paper \cite{Ma-Pa} for a detailed list of examples.
\end{remark}


\subsection{Hessian comparison analysis}\label{subsecLap}
The 2.nd order analysis of the restricted distance function $r_{|_{P}}$ defined on manifolds with a pole is  governed by the Hessian comparison Theorem A in \cite{GreW}.

This comparison theorem can be stated as follows,  when one of the spaces is a model space $M^m_w$, (see \cite{Pa3}):

\begin{theorem}[See \cite{GreW}, Theorem A]\label{thmGreW}
Let $N=N^{n}$ be a manifold with a pole $o$, let $M=M_{w}^{m}$ denote a
$w-$model with center $o_{w}$. Suppose that every $o$-radial
sectional curvature at $x \in N \setminus \{o\}$ is bounded from above by
the $o_{w}$-radial sectional curvatures in $M_{w}^{m}$ as follows:
\begin{equation*}
K_{o, N}(\sigma_{x})\, \leq\, -\frac{w''(r)}{w(r)}
\end{equation*}
for every radial two-plane $\sigma_{x} \in T_{x}N$ at distance $r =
r(x) = \dist_{N}(o, x)$ from $o$ in $N$. Then the Hessian of the
distance function in $N$ satisfies
\begin{equation}\label{eqHess}
\begin{aligned}
\Hess^{N}(r(x))(X, X) &\, \geq\, \Hess^{M}(r(y))(Y, Y)\\ &=
\eta_{w}(r)\left(\Vert X\Vert^2 - \langle \nabla^{M}r(y), Y \rangle_{M}^{2}
\right) \\ &= \eta_{w}(r)\left(\Vert X\Vert^2 - \langle \nabla^{N}r(x), X
\rangle_{N}^{2} \right)
\end{aligned}
\end{equation}
for every  vector $X$ in $T_{x}N$ and for every vector $Y$
in $T_{y}M$ with $\,r(y) = r(x) = r\,$ and $\, \langle
\nabla^{M}r(y), Y \rangle_{M} = \langle \nabla^{N}r(x), X
\rangle_{N}\,$.
\end{theorem}

\begin{remark}
As we mentioned in the Introduction, inequality (\ref{eqHess}) is true along the geodesics emanating from $o$ and $o_w$ which are free of conjugate points of $o$ and $o_w$, (see Remark 2.3 in \cite{GreW}). Other relevant observation is that the bound given in inequality (\ref{eqHess}) does not depend on the dimension of the model space, (see Remark 3.7 in \cite{Pa3}). 
\end{remark}

We present now a technical result concerning the Hessian of a radial function, namely, a function which only depends on the distance function $r$. For the proof of this result, and the rest of the results in this subsection, we refer to the paper \cite{Pa3}.
 \medskip
 
 \begin{proposition}\label{hessianfunction}
 Let $N=N^{n}$ be a manifold with a pole $o$. Let $r =
r(x) = \dist_{N}(o, x)$ be the distance from $o$ to $x$ in $N$. Let $F:\erre \longrightarrow \erre$ a smooth function. Then, given $q\in N$ and  $X, Y \in T_qN$,
\begin{equation}\label{hessF}
\begin{split}
\Hess^N F\circ r\vert_q(X,Y)&=F''(r)(\nabla^N r \otimes \nabla^N r)(X,Y)\\&+F'(r)\Hess^N r\vert_q(X,Y)
\end{split}
\end{equation}
 \end{proposition}

Now, let us consider a complete isometric immersion $\varphi: P^m \longrightarrow N$ in a Riemannian ambient manifold $N^n$ with pole $o$, and with distance function to the pole $r$.  We are going to see how the Hessians (in $P$ and in $N$),  of a radial function defined in the submanifold are related via the second fundamental form $B^P$ of the submanifold $P$ in $N$.  As before, we identify, given any $q \in P$, the tangent vectors $X \in T_qP$ with $\varphi_{*_q}X \in T\varphi(q)N$ along the next results.

\begin{proposition}
Let $N^{n}$ be a manifold with a pole $o$, and let us consider an isometric immersion $\varphi: P^m \longrightarrow N$. If $r\vert_P$ is the extrinsic distance function, then, given $q\in P$ and $X, Y \in T_qP$,
\begin{equation}
\Hess^P r\vert_q(X,Y)=\Hess^N r\vert_{\varphi(q)}(X,Y) +\langle B^P_{q}(X,Y),\nabla^N r\vert_q\rangle
\end{equation}
where $B^P_{q}$ is the second fundamental form of $P$ in $N$ at the point $q \in P$.
\end{proposition}

Now, we apply Proposition \ref{hessianfunction} to $F\circ r\vert_P=F\circ r\circ \varphi$, (considering $P$ as the Riemannian manifold where the function is defined), to obtain an expression for $\Hess^P F\circ r\vert_P(X,Y)$ . Then, let us apply Proposition above to $\Hess^P r\vert_P(X,Y)$, and we finally get:

\begin{proposition}\label{hessianfunctionrestricted}
Let $N=N^{n}$ be a manifold with a pole $o$, and let $P^m$ denote an immersed submanifold in $N$. Let $r\vert_P$ be the extrinsic distance function. Let $F:\erre \longrightarrow \erre$ be a smooth function. Then, given $q\in P$ and  $X, Y \in T_qP$,
\begin{equation}\label{subhessf}
\begin{split}
\Hess^P F\circ r\vert_q(X,Y)&=F''(r(q))\langle\, \nabla^N r\vert_q,X\, \rangle\langle\, \nabla^N r\vert_q,Y\, \rangle\\&+F'(r(q))\{\Hess^N r\vert_q(X,Y)\\&+\langle\nabla^N r\vert_q,B^P_q(X,Y)\, \rangle \, \}
\end{split}
\end{equation}
 \end{proposition}


\subsection{Comparison constellations and Isoperimetric inequalities}

The isoperimetric inequalities satisfied by the extrinsic balls in minimal submanifolds are on the basis of the monotonicity of the volume growth function $f(r)=\frac{Vol(D_r)}{Vol(B_r^w)}$, a key result to prove Theorem \ref{cor1}. We have the following theorem.


\begin{theorem}[See \cite{Ma-Pa}, \cite{MP}, \cite{MP1}, \cite{MP2} and \cite{Pa1}]\label{Isoperimetric}

Let $\varphi: P^m \longrightarrow N^n$ be a complete, proper and minimal immersion in an ambient Riemannian manifold $N^n$ which possess at least one pole $o \in N$.  Let us suppose that the $o-$radial sectional curvatures of $N$ are bounded from above by the $o_w-$radial sectional curvatures of the $w-$model space $M_w^m$:
$$
K_{o,N}(\sigma_{x}) \leq  -\frac{w''(r(x))}{w(r(x))}\,\,\,\forall x \in N
$$

\noindent  and assume that $M^{m}_{w}$ is balanced from below.
Let $D_r$ be an extrinsic $r$-ball in $P^m$,
 with center at a pole $ o \in N$ in the ambient space $N$. Then:
 
\begin{equation} \label{isopComp}
\frac{\Vol(\partial D_r)}{\Vol(D_r)} \geq
\frac{\Vol(S^{w}_r)}{\Vol(B^{w}_r)} \,\,\,\,\,\textrm{for all}\,\,\,
r>0  \quad .
\end{equation}
Furthermore, if $\varphi^{-1}(o)\neq \emptyset$,
    \begin{equation} \label{eqVolComp}
    \Vol(D_r) \geq \Vol(B^w_r)\,\,\,\,\textrm{for all}\,\,\, r>0
    \quad .
    \end{equation}
Moreover, if equality in inequalities (\ref{isopComp}) or
(\ref{eqVolComp}) holds
    for some fixed radius $R$ and if 
    the balance of $M^{m}_{w}$ from below is sharp $q_{w}(r)\,\eta_{w}(r)\, > \, 1/m\, $ for all $r$,
then $D_R$ is  a
    minimal cone in the ambient space $N^n$, so if $N^n$ is the
    hyperbolic space $\,\Han\,$, $\,b < 0\,$, then $P^m\,$ is totally geodesic in
    $\Han$.
    
      If, on the other hand, the ambient space is $\erre^n$ and equality in inequalities (\ref{isopComp}) or
(\ref{eqVolComp}) holds
    for all radius $r>0$ 
then $P^m$ is totally geodesic in
    $\erre^n$.

    \medskip
    
    On the other hand, and also as a consequence of inequality (\ref{isopComp}), the volume growth function $f(r)=\frac{Vol(D_r)}{Vol(B_r^w)}$ is a non-decreasing function of $r$.
\end{theorem}


\section{Main Results}\label{Main-Results}
 We prove in this section our main results, stablishing a set of conditions that assures that our submanifolds are properly immersed and have finite topology and bounding from below, under certain conditions, the number of its ends. 

\begin{theorem}\label{theorem1} 
Let $\varphi: P^m \longrightarrow N^n$ be an isometric immersion of a complete non-compact Riemannian $m$-manifold $P^m$ into a complete Riemannian manifold $N^n$ with a pole $o\in N$ and satisfying $\varphi^{-1}(o) \neq \emptyset$. Let us suppose that:
\begin{enumerate}
\item The $o-$radial sectional curvatures of $N$ are bounded from above by the $o_w-$radial sectional curvatures of the $w-$model space $M_w^m$:
$$
K_{o,N}(\sigma_{x}) \leq  -\frac{w''(r(x))}{w(r(x))}\,\,\,\forall x \in N.
$$
\item The second fundamental form $B^P_x$ in $x\in P$ satisfies that, for sufficiently large radius $R_0$, and for some constant $c \in ]0,1[$:
$$\Vert B^P_x \Vert\leq c\, \eta_w(\rho^P(x))\,\,\,\forall x\in P - B^P_{R_0}(x_o)$$
where $\rho^P(x)$ denotes the intrinsic distance in $P$ from some fixed $x_o \in \varphi^{-1}(o)$ to $x$.
\item For any $r>0$,  $w'(r)\geq d>0$ and $(\eta_w(r))'\leq 0$.
\end{enumerate}

Then $P$ is properly immersed in $N$ and it is $C^{\infty}$-
diffeomorphic to the interior of a  compact smooth manifold $\overline{P}$ with
boundary.
\end{theorem}

\begin{remark} To show that $\varphi$  is proper, we shall use Theorem \ref{thmGreW}. Hence, it is enough to assume that $o$ is a pole in the sense that there are not conjugate points along any geodesic emanating from $o$, (see \cite{doCarmo} and \cite{S}). Therefore our statement about the properness of the immersion  includes ambient manifolds $N$ that admit non-negative sectional curvatures, unlike the ambient manifold in Theorem 1.2 in \cite{Pac2}. On the other hand, to prove the finiteness of the topology of $P$ we need to assume that the ambient manifold $N$ posses a pole as it is defined in \cite{GreW}, namely, a point $p \in N$ where $exp_p$ is a $C^{\infty}$ diffeomorphism. However, although our ambient manifold must be diffeomorphic to $\erre^n$ in this case, (as in Theorem 1.2 in \cite{Pac2}, where the ambient space must be a Cartan-Hadamard manifold), also admits non-negative sectional curvatures.

To complete the benchmarking with the hypotheses in \cite{Pac} and \cite{Pac2}, we are going to compare the assumptions (2) and (3) in Theorem \ref{theorem1} with the notion of \lq\lq submanifold with {\it tamed} second fundamental form" introduced in \cite{Pac}. It is straightforward to check that if $\varphi: P^m \longrightarrow N^n$ is an immersion of a complete Riemannian $m$- manifold $P^m$ into a complete Riemannian manifold $N^n$  with sectional curvatures $K_N\leq b\leq 0$, and $P$ has tamed second fundamental form, in the sense of Definition 1.1 in \cite{Pac2}, then there exists $R_0>0$ such that for all $r \geq R_0$, the quantity 
$$a_r:= \Sup\{ \frac{w_b}{w_b'}(\rho^P(x))\Vert B^P_x\Vert: x \in P -B^P_{r}\}$$ 
\noindent satisfies $a_r <1$.

Hence, taking $r=R_0$, we have that for all $x \in P -B^P_{R_0}$, and some $c \in (0,1)$,

$$\Vert B^P_x\Vert \leq c \eta_{w_b}(\rho^P(x))\,.$$

On the other hand, when $b\leq 0$, then $w_b'(r) \geq 1 >0\,\,\forall r>0$ and $(\eta_{w_b}(r))'\leq 0\,\,\forall r>0$.

All these observations make us consider our Theorem \ref{theorem1} as a natural and slight generalization of assertions (b) and (c) of Theorem 1.2 in \cite{Pac2}.

\end{remark}

Observe that if we assume the properness of the immersion we obtain the following version of Theorem \ref{theorem1}, where we can remove the hypothesis about the decrease of the function $\eta_w(r)$ because the norm of the second fundamental form $\Vert B^P_x \Vert$ is bounded by the value of $\eta_w$ at $r(x)$ instead of $\rho^P(x)$ :
\begin{theorem}\label{alternative}
Let $\varphi: P^m \longrightarrow N^n$ be an isometric and proper immersion of a complete non-compact Riemannian $m$-manifold $P^m$ into a complete Riemannian manifold $N^n$ with a pole $o\in N$ and satisfying $\varphi^{-1}(o) \neq \emptyset$.  Let us suppose that, as in Theorem \ref{theorem1}, the $o-$radial sectional curvatures of $N$ are bounded from above as 
$$K_{o,N}(\sigma_{x}) \leq  -\frac{w''(r(x))}{w(r(x))}\,\,\,\forall x \in N\, ,$$ and for any $r>0$,  $w'(r)\geq d>0$. Let us assume moreover that the second fundamental form $B^P_x$ in $x\in P$ satisfies that, for sufficiently large radius $R_0$: 
$$\Vert B^P_x \Vert\leq c\, \eta_w(r(x))\,\,\,\forall x\in P - D_{R_0}(o)$$
where $c$ a positive constant such that $c<1$ .

Then $P$  is $C^{\infty}$-
diffeomorphic to the interior of a  compact smooth manifold $\overline{P}$ with
boundary.
\end{theorem}

We are going to see how to estimate the area growth function of $P$, defined as $g(r)=\frac{Vol(\partial D_r)}{Vol(S_r^w)}$ by the number of ends of the immersion $P$, $\mathcal{E}(P)$, when the ambient space $N$ is a radially symmetric space.

\begin{theorem}\label{theorem2}  Let $\varphi: P^m \longrightarrow M^n_w$ be an isometric and proper immersion of a complete non-compact Riemannian $m$-manifold $P^m$ into a model space $M^n_w$ with pole $o_w$. Suppose that $\varphi^{-1}(o_w)\neq \emptyset$, $m>2$ and moreover:
  \begin{enumerate}

\item The norm of second fundamental form $B^P_x$ in $x\in P$ is bounded from above outside a (compact) extrinsic ball $D_{R_0}(o)\subseteq P$ with sufficiently large radius $R_0$  by: 

$$\Vert B^P_x \Vert\, \leq\, \frac{\epsilon(r(x))}{(w'(r(x)))^2}\eta_w(r(x))\,\,\,\forall x\in P - D_{R_0}$$
\noindent where  $\epsilon$ is a positive function such that $\epsilon(r)\to 0$ when $r\to \infty$. 

\item For $r$ sufficiently large, $w'(r)\geq d>0$.
\end{enumerate}

 \noindent Then,  for sufficiently large $r$, we have:
 \begin{equation}\label{lowerboundends}
 \frac{Vol(\partial D_r)}{Vol(S_r^w)}\leq \frac{\mathcal{E}(P)}{\left(1-4\epsilon(r)\right)^{\frac{(m-1)}{2}}}
 \end{equation}

\noindent where $\mathcal{E}(P)$ is the (finite) number of ends of $P$.

\end{theorem}

When we consider minimal immersions in the model spaces, we have the following result, which is an inmediate corollary from the above theorem,  and Theorem \ref{Isoperimetric}  in Section \ref{Prelim}.
\newpage

\begin{theorem}\label{minimal}
 Let $\varphi: P^m \longrightarrow M^n_w$ be a complete non-compact, proper and minimal immersion into a ballanced from below model space $M^n_w$ with pole $o_w$. Suppose that $\varphi^{-1}(o_w)\neq \emptyset$ and  $m>2$.  Let us assume moreover the hypotheses (1) and (2) in Theorem \ref{theorem2}.
 
Then 
\begin{enumerate}
\item The (finite) number of ends $\mathcal{E}(P)$ is related with the (finite) volume growth by
\begin{equation}\label{volgrowthends}
1\leq \lim_{r\to \infty}\frac{Vol(D_r)}{Vol(B_r^w)}\leq\mathcal{E}(P)
 \end{equation}

\item If $P$ has only one end, P is a minimal cone  in $M_w^n$.
\end{enumerate}
\end{theorem}


\section{Corollaries}\label{corollaries}

As we have said in the Introduction, we have divided the list of results based in Theorem \ref{theorem1} and in Theorem \ref{theorem2} in two series of corollaries. The first set of consequences follows the line of Theorem \ref{cor1} and Theorem \ref{cor3}, (which are in fact the main representatives of these results) presenting upper bounds for the volume and area growth of a complete and proper immersion in the real space form $\kan$, ($b \leq 0$), in terms of the number of its ends. In the second set of corollaries, are stated compactification theorems for complete and proper immersions in $\erre^n$, $\Han$ and $\Han \times \erre^l$.

The first of these corollaries constitutes a non-minimal version of Theorem \ref{cor1}:

\begin{corollary}\label{cor2}
Let $\varphi: P^m \longrightarrow \Han$ be a complete non-compact and proper immersion with $m>2$. Let us suppose that for sufficiently large $R_0$ and for all points $x \in P$ such that $r(x)>R_0$, (i.e. outside the compact extrinsic ball $D_{R_0}(o)$ with $\varphi^{-1}(o) \neq \emptyset$),

 $$\Vert B^P_x \Vert \leq \frac{\delta(r(x))}{e^{2\sqrt{-b}\, r(x)}}$$
 
\noindent where $r(x)=d_{\mathbb{H}^n(b)}(o,\varphi(x))$ is the (extrinsic) distance in $\mathbb{H}^n(b)$ of the points in $\varphi(P)$ to a fixed pole $o \in \mathbb{H}^n(b)$ and $\delta(r)$ is a smooth function such that $\delta(r)\to 0$ when $r\to \infty$.Let $\{t_i\}_{i=1}^\infty$ be any non-decreasing sequence such that $t_i\to \infty$ when $i\to \infty$. Then the finite number of ends $\mathcal{E}(P)$ is related with the area growth of $P$ by:
$$\liminf_{i \to \infty}\frac{\Vol(\partial D_{t_i})}{\Vol(S_{t_i}^{b,m-1})}\leq \mathcal{E}(P)$$
\end{corollary}

The corresponding non-minimal statement of Theorem \ref{cor3} is:

\begin{corollary}\label{cor4}
Let $\varphi: P^m \longrightarrow \erre^n$ be a complete non-compact and proper immersion with $m>2$. Let us suppose that for sufficiently large $R_0$ and for all points $x \in P$ such that $r(x)>R_0$, (i.e. outside the compact extrinsic ball $D_{R_0}(o)$ with $\varphi^{-1}(o) \neq \emptyset$),
 $$\Vert B^P_x \Vert \leq \frac{\epsilon(r(x))}{r(x)}$$
where $r(x)=d_{\erre^n}(o,\varphi(x))$ is the (extrinsic) distance in $\erre^n$ of the points in $\varphi(P)$ to a fixed pole $o \in \erre^n$ and $\epsilon(r)$ is a smooth function such that $\epsilon(r)\to 0$ when $r\to \infty$. Let $\{t_i\}_{i=1}^\infty$ be any non-decreasing sequence such that $t_i\to \infty$ when $i\to \infty$. Then the finite number of ends $\mathcal{E}(P)$ is related with the area growth by:
$$\liminf_{i \to \infty}\frac{\Vol(\partial D_{t_i})}{\Vol(S_{t_i}^{0,m-1})}\leq \mathcal{E}(P)$$
\end{corollary}

Concerning the compactification results we have the following result given by Bessa, Jorge and Montenegro in \cite{Pac} and by Bessa and Costa in \cite{Pac2}:

\begin{corollary}\label{cor5}
Let $\varphi: P^m \longrightarrow \kan$ be a complete non-compact immersion in the real space form $\kan$, ($b \leq 0$). Let us suppose that for all points $x\in P\setminus B^P_{R_0}(o)$ (for sufficientlty large $R_0$, where $o$ is a pole in $\kan$ such that $\varphi^{-1}(o) \neq \emptyset $) :
 $$\Vert B^P_x \Vert \leq c\,h_b(\rho^P(x))$$
where $\rho^P(x)$ is the (intrinsic) distance to a fixed $x_o \in \varphi^{-1}(o)$ and  $c$ is a positive constant such that $c<1$ and 
$$h_b(r)= \eta_{w_b}(r)=\begin{cases} 
\phantom{\sqrt{b}} 1/r &\text{if $b=0$}\\\sqrt{-b}\coth(\sqrt{-b}\,r) &\text{if
$b<0$} \quad .
\end{cases}$$
is the mean curvature of the geodesic spheres in $\kan$.
Then $P$ is properly immersed in $\kan$ and it is diffeomorphic to the interior of a compact smooth manifold $\overline P$ with boundary. 
\end{corollary}

Our last result concerns isometric immersions in $\mathbb{H}^n(b)\times \erre^ l$:

\begin{corollary}\label{cor6}
Let $\varphi: P^m \longrightarrow \mathbb{H}^n(b)\times \erre^ l$ be a complete non-compact immersion. Let us consider a pole $o \in \mathbb{H}^n(b)\times \erre^ l$ such that $\varphi^{-1}(o) \neq \emptyset$. Let us suppose  that for all points $x\in P\setminus B^P_{R_{0}}(x_o)$, where $x_o \in \varphi^{-1}(o)$ and for $R_0$ sufficiently large:
$$\Vert B_x \Vert \leq \frac{c}{\rho^P(x)}\,\, .$$
\noindent Here $\rho^P(x)$ denotes the intrinsic distance in $P$ from the fixed $x_o \in \varphi^{-1}(o)$ to $x$ and $c$ is a positive constant such that $c<1$. Then $P$ is properly immersed in $\mathbb{H}^n(b)\times \erre^ l$ and it is diffeomorphic to the interior of a compact smooth manifold $\overline P$ with boundary. 
\end{corollary}


\section{Proof of Theorem \ref{theorem1}}\label{demos1}


\subsection{$P$ is properly immersed}

 Let us define the following function:
\begin{equation}
F(r):=\int_0^r w(t)dt
\end{equation}
Observe that $F$ is injective, because $F'(r)=w(r)> 0\,\,\,\forall r>0$, and  $F(r)\to \infty$ when $r\to \infty$. Applying Theorem \ref{thmGreW} and Proposition \ref{hessianfunctionrestricted}, we obtain, for all $x \in P$, and given $X \in T_xP$,

\begin{equation}\label{convex-hess1}
\begin{aligned}
\Hess^P_x F(r)(X,X) &\geq   w'(r(x))\Vert X \Vert^2+w(r(x))\langle B^P_x(X,X),\nabla^N r\rangle \\&\geq w'(r(x))\Vert X \Vert^2-w(r(x))\Vert B^P_x\Vert\,\, \Vert X\Vert ^2
\end{aligned}
\end{equation}

\noindent By hypotesis there exist a geodesic ball $B^P_{r_1}(x_0)$ in $P$, with $r_1 \geq R_0$, such that for any $x\in P\setminus B^P_{r_1}(x_0)$, $\Vert B^P_x\Vert\ \leq c\eta_w(\rho^P(x))$. On the other hand, as $\eta_w(r)$ is non-increasing and $r(x) \leq \rho^P(x)$ because $\varphi$ is isometric, we have $c\eta_w(\rho^P(x))\leq c\eta_w(r(x))$, so if $x\in P\setminus B^P_{r_1}$ :

\begin{equation}\label{convex-hess}
\begin{aligned}
\Hess^P_x F(r)(X,X) &\geq 
w'(r(x))\Vert X \Vert^2-w(r)c\eta_w(\rho^P(x))\, \Vert X\Vert ^2\\&\geq w'(r(x))\Vert X \Vert^2\left(1-c\right)\geq d\left(1-c\right)>0
\end{aligned}
\end{equation}

\noindent The above result implies that there exists $r_1 \geq R_0$ such that $F\circ r$ is a strictly convex function  outside the geodesic ball in $P$ centered at $x_0$, $B^P_{r_1}(x_0)$. And hence, as $r(x) \leq \rho^P(x)$ for all $x \in P$, (and therefore $B^P_{r_1}(x_0) \subseteq D_{r_1}$), $F\circ r$ is a strictly convex function outside the extrinsic disc $D_{r_1}$. 

 Let $\sigma :[0,\rho^P (x)]\to P^{m}$ be a minimizying
geodesic from $x_0$ to $x$. 

If we denote as $f=F\circ r$, let us define $h: \erre \to \erre$ as 
$$h(s)=F(r(\sigma(s)))=f(\sigma(s))$$

 Then, 
\begin{equation}
(f \circ \sigma)'(s)=h'(s)=\sigma'(s)(f)=\langle \nabla^P f(\sigma(s)),\sigma'(s)\rangle
\end{equation}
and hence,
\begin{equation}
\begin{aligned}
(f \circ \sigma)''(s)&=h''(s)=\sigma'(s)(\langle\nabla^P f(\sigma(s)),\sigma'(s)\rangle)=\langle\nabla^P_{\sigma'(s)}\nabla^P f(\sigma(s)),\sigma'(s)\rangle\\&+\langle\nabla^P f(\sigma(s)),\nabla^P_{\sigma'(s)}\sigma'(s)\rangle=Hess^P_{\sigma(s)}f(\sigma(s))(\sigma'(s),\sigma'(s))
\end{aligned}
\end{equation}

We have from (\ref{convex-hess}) that 
 $(f\circ\sigma)''(\tau)=\Hess^P f(\sigma
(\tau))(\sigma', \sigma')\geq d(1-c)$ for all
$\tau\geq r_1$ . And for $\tau<r_1$,
$(f\circ\sigma)''(\tau))\geq a=\inf_{x\in B^P_{r_1}} \{ \Hess^P
f(x)(\nu,\nu), \,\vert \nu\vert=1\}$. Then
\begin{eqnarray}\label{eq2}
(f\circ\sigma) ' (s)&=&(f\circ\sigma) ' (0)+ \int_{0}^{s}(f \circ \sigma)'' (\tau)d\tau \nonumber  \\
&\geq & (f\circ\sigma) ' (0)+\int_{0}^{r_1}a\,d\tau +d\, \int_{r_1}^{s} (1-c)d\tau\\
&&\nonumber \\ &\geq &(f\circ\sigma) ' (0)+ a\,r_1+ d\,(1-c) (s-r_1)\nonumber
\end{eqnarray}

\noindent On the other hand, as 
\begin{equation}
\nabla^P f(\sigma(s))=\nabla^P F(r(\sigma(s)))=F'(r(\sigma(s)))\nabla^Pr\vert_{\sigma(s)}=w(r(\sigma(s)))\nabla^Pr\vert_{\sigma(s)}
\end{equation}

\noindent then 
$$\nabla^P f(\sigma(0))=w(r(\sigma(0)))\nabla^Pr\vert_{\sigma(0)}=w(0)\nabla^Pr\vert_{\sigma(0)}=0$$
\noindent so we have that
\begin{equation}
(f\circ\sigma)'(0)=\langle \nabla^P f(\sigma(0)),\sigma'(0)\rangle=0
\end{equation}

\noindent We also have that $(f\circ\sigma)(0)=F(r(\sigma(0)))=F(0)=0$.
Hence, applying inequality (\ref{eq2}),

\begin{equation}
f(\sigma(s))=(f\circ \sigma)(0)+ \int_0^s (f\circ \sigma)'(\tau)d\tau \geq ar_1s+d(1-c)\{\frac{1}{2} s^2-r_1s\}
\end{equation}

Therefore, 
\begin{eqnarray}
F(r(x))&=&f(x)=f(\sigma(\rho^P(x)))= \int_{0}^{\rho^P(x)}(f\circ \sigma)' (s) \, ds\nonumber \\
 &\geq & \int_{0}^{\rho^P(x)} a\,r_1+ d\,(1-c) (s-r_1)\,ds \label{eq3} \\
 &=&a\,r_1 \rho^P(x)\ + d\,(1-c)\left( \frac{\rho^P(x)^{2}}{2}-r_1\,\rho^P(x)\right)\nonumber
\end{eqnarray} 

Hence, if $\rho^P\to \infty$ then $F(r(x))\to \infty$ and then, as  $F$ is strictly increasing, $r\to \infty$  so the immersion is proper. 


\subsection{$P$ has finite topology}

We are going to see that $\nabla^ Pr$ never vanishes on 
$P\setminus D_{r_1}$. To show this, we consider, as in the previous subsection, any geodesic in $P$ emanating from the pole $o$, $\sigma(s)$. We have, using inequality (\ref{eq2}), that 
\begin{equation}
\langle \nabla^P f(\sigma(s)),\sigma'(s)\rangle= (f \circ \sigma)'(s) \geq a\,r_1+ d\,(1-c) (s-r_1) >0\,\,\forall s>r_1
\end{equation}

Hence, as $\Vert \sigma'(s)\Vert=1 \,\,\forall s$, then $\Vert \nabla^P f(\sigma(s))\Vert >0$ for all $s >r_1$. But we have computed $\nabla^P f(\sigma(s))=w(r(\sigma(s)))\nabla^Pr\vert_{\sigma(s)}$, so, as $w(r)>0\,\,\forall r>0$, then $\Vert \nabla^Pr\vert_{\sigma(s)}\Vert >0\,\,\forall s>r_1$ and hence, $\nabla^Pr\vert_{\sigma(s)}\neq 0\,\,\forall  s>r_1$. We have proved that $\nabla^ Pr$ never vanishes on 
$P\setminus B^P_{r_1}$, so we have too that  $\nabla^ Pr$ never vanishes on 
$P\setminus D_{r_1}$. 
Let 
$$
\phi : \partial D_{r_1}\times [r_1,+\infty)\to P\setminus D_{r_1}
$$
be the integral flow of a vector field $\frac{\nabla^P r}{\Vert \nabla^P r \Vert ^2}$ with
$$
\phi(p,r_1)=p\in \partial D_{r_1}
$$ 
It is obvious that $r(\phi(p,t))=t$ and
$$
\phi(\cdot ,t):\partial D_{r_1} \to \partial D_{t}
$$
is a diffeomorphism. So $P$ has finitely many ends, and each of its ends is of finite topological type.

In fact, applying Theorem 3.1 in \cite{Mi}, we conclude that, as the extrinsic annuli $A_{r_1,R}(o)= D_R(o) \setminus D_{r_1}(o)$ contain no critical points of the extrinsic distance function $r: P \longrightarrow \erre^+$, then $D_R(o)$ is diffeomorphic to $D_{r_1}(o)$ for all $R \geq r_1$ and hence the annuli $A_{r_1,R}(o)$ are diffeomorphic to $\partial D_{r_1} \times [r_1, R]$.

\begin{remark}
To show Theorem \ref{alternative}, we argue as in the beginning of the proof of Theorem \ref{theorem1}: with the same function $F(r)$ we obtain inequality (\ref{convex-hess1}). But now we have as hypothesis that $\Vert B^P_x \Vert\leq c\, \eta_w(r(x))$, so we don't need that $\eta_w'(r) \leq 0$ to get inequality (\ref{convex-hess}).
\end{remark}


\section{Proof of Theorem \ref{theorem2}}\label{demos2}

We are going to see first that $P$ has finite topology. As $P$ is properly immersed, we shall apply Theorem \ref{alternative} and for that, it must be checked that hypotheses in that theorem are acomplished. First, we have hypothesis (1) in Theorem \ref{alternative} because $N=M^n_w$. On the other hand, as $w'(r) \geq d >0 \forall r >0$ and, for some $R_0$, we have that $\Vert B^P_x\Vert \leq \frac{\epsilon(r(x))}{(w'(r(x)))^2}\eta_w(r(x))\,\,\,\forall x\in P - D_{R_0}$ where  $\epsilon$ is a positive function such that $\epsilon(r)\to 0$ when $r\to \infty$, hence $0 \leq \lim_{r \to \infty} \frac{\epsilon(r)}{(w'(r))^2} \leq \lim_{r \to \infty} \frac{\epsilon(r)}{d^2}=0$. Therefore, for some constant $c <1$, there exist $R_0$ such that $\Vert B^P_x\Vert \leq c\eta_w(r(x))\,\,\,\forall x\in P - D_{R_0}$.
Therefore, as $\varphi: P \longrightarrow M^n_w$ is a proper immersion, we have by Theorem \ref{alternative}  that  $P$ has finite topological type and thus $P$ has finitely many ends, each of finite topological type. Hence we have, in an analogous way than in \cite{anderson}, and for $r_1 \geq R_0$ as in Section \ref{demos1}:

\begin{equation}
P - D_{r_1}=\cup_{k=1}^{\mathcal{E}(P)} V_k
\end{equation}
where $V_k$ are disjoint, smooth domains in $P$. Along  the rest of the proof, we will work on each end $V_k$ separately. Let $V$ denote one element of the family  $\{ V_k\}_{k=1}^{\mathcal{E}(P)}$, and, given a fixed radius $t > r_1$, let $\partial V(t)$ denote the set $\partial V(t)=V\cap \partial D_t=V\cap S^w_t$, where $S^w_t$ is the geodesic $t$-sphere in $M^n_w$. This set is a hypersurface in $P^m$, with normal vector $\frac{\nabla^Pr}{\Vert \nabla^P r\Vert}$, and we are going to estimate its sectional curvatures when $t \to \infty$.

 
  Suppose that $e_i,e_j$ are two orthonormal vectors of $T_p\partial V(t)$ on the point $p\in \partial V(t)$. Then the sectional curvature of the plane expanded by $e_i,e_j$ is, using Gauss formula:
  \begin{equation}
  \begin{aligned}
  K_{\partial V(t)}&(e_i,e_j)=K_P(e_i,e_j)+\langle B^{\partial V-P}(e_i,e_i),B^{\partial V-P}(e_j,e_j)\rangle\\&-\Vert B^{\partial V-P}(e_i,e_j)\Vert ^2=K_N(e_i,e_j)+\langle B^{\partial V-P}(e_i,e_i),B^{\partial V-P}(e_j,e_j)\rangle\\&-\Vert B^{\partial V-P}(e_i,e_j)\Vert ^2+\langle B^{P}(e_i,e_i),B^{P}(e_j,e_j)\rangle-\Vert B^{P}(e_i,e_j)\Vert ^2\\&\geq K_N(e_i,e_j)+ \langle B^{\partial V-P}(e_i,e_i),B^{\partial V-P}(e_j,e_j)\rangle\\&-\Vert B^{\partial V-P}(e_i,e_j)\Vert ^2-2\Vert B^{P}\Vert ^2
  \end{aligned}
  \end{equation}
 
\noindent where $B^{\partial V-P}$ is the second fundamental form of $\partial V(t)$ in $P$. But this second fundamental form is for two vector fields $X,Y$ in $T\partial V(t)$:
\begin{equation}
\begin{aligned}
B^{\partial V-P}(X,Y)&=\langle\nabla_X^PY,\frac{\nabla ^Pr}{||\nabla ^Pr||}\rangle\frac{\nabla ^Pr}{||\nabla ^Pr||}=\langle\nabla_X^PY,\nabla ^Pr\rangle\frac{\nabla ^Pr}{||\nabla ^Pr||^2}\\&=X(\langle Y,\nabla^Pr\rangle)\frac{\nabla ^Pr}{||\nabla ^Pr||^2}-\langle Y,\nabla^P_X\nabla^Pr\rangle\frac{\nabla ^Pr}{||\nabla ^Pr||^2}\\&=-\Hess^Pr(X,Y)\frac{\nabla ^Pr}{||\nabla ^Pr||^2}
\end{aligned}
\end{equation}
Then, since, for all $X,Y \in T_pM^n_w$
\begin{equation}
\Hess^{M^{n}_w} r(X,Y)=\eta_w(r)\langle X,Y \rangle- \langle X, \nabla^{M^{n}_w} r \rangle\langle Y, \nabla^{M^{n}_w} r \rangle
\end{equation}

\noindent we have, (using the fact that $e_i$ are tangent to the fiber $S^w_t$, and Proposition \ref{propWarpMean}), that 
\begin{equation}
K_{M^n_w}(e_i,e_j)=K(t)=\frac{1}{w^2(t)}-\eta^2_w(t)
\end{equation}

so for  any $p\in \partial V(t)$ such that $t=r(p)$ is sufficiently large:

\begin{equation}
  \begin{aligned}
  K_{\partial V(t)}(e_i,e_j)\geq&K_{M^n_w}(e_i,e_j)+ \frac{\Hess^P_pr(e_i,e_i)\Hess^P_pr(e_j,e_j)}{||\nabla ^Pr||^2}\\&-\frac{\Hess^P_pr(e_i,e_j)^2}{||\nabla ^Pr||^2}-2\Vert B^{P}\Vert ^2 
  \\&\geq K(t)+\frac{\left(\eta_w(t)-\Vert B^ P \Vert\right)^ 2-\Vert B^P \Vert ^2}{||\nabla ^Pr||^2}-2\Vert B^{P}\Vert ^2\\&
  \geq \eta^2_w(t)\left(1-2\frac{\Vert B^ P \Vert}{\eta_w(t)}-2\left(\frac{\Vert B^ P \Vert}{\eta_w(t)}\right) ^2+\frac{K(t)}{\eta^2_w(t)}\right)\\&
  \geq \eta^2_w(t)\left(1-4\frac{\Vert B^ P \Vert}{\eta_w(t)}+\frac{K(t)}{\eta^2_w(t)}\right)\\&
  =\eta^2_w(t)\left( 1+\frac{K(t)}{\eta^2_w(t)}\right)\left(1-4\frac{\frac{\Vert B^ P \Vert}{\eta_w(t)}}{1+\frac{K(t)}{\eta^2_w(t)}}\right)\\&
  \geq\frac{1}{w^2(t)}\left(1-4\Vert B^ P \Vert w'(t)w(t)\right)\geq\frac{1}{w^2(t)}\left(1-4\epsilon(t)\right)
  \end{aligned}
  \end{equation}
  \noindent where we recall that, by hypothesis, $\Vert B^ P \Vert \leq \frac{\epsilon(t)}{(w'(t))^2}\eta_w(t)$ for all $t=r(x)>R_0$, and $\epsilon$ is a positive function such that $\epsilon(r)\to 0$ when $r\to \infty$.
  
  If we  denote as $\delta(t)=\frac{1}{w^2(t)}\left(1-4\epsilon(t)\right)$ we 
 have for each $t$ sufficiently large that $K_{\partial V(t)}(e_i,e_j)\geq \delta(t)$ holds everywhere on $\partial V(t)$ and $\delta(t)$ is a positive constant. Then, the Ricci curvature of $\partial V(t)$ is bounded from below, for these sufficiently large radius $t$ as
 $$Ricc_{\partial V(t)}(\xi,\xi)\geq \delta(t)(m-1)\Vert \xi\Vert^2 >0\,\,\forall \xi \in T\partial V(t)$$
 
\noindent so, applying Myers' Theorem  $\partial V(t)$ is compact and has diameter $d(\partial V(t))\leq \frac{\pi}{\sqrt{\delta(t)}}$ (see \cite{S}). Applying on the other hand Bishop's Theorem, (see Theorem 6 in \cite{Chavel}), we obtain:
 
 \begin{equation}
  \begin{aligned}
 \Vol(\partial V(t))\leq \frac{\Vol(S^{0,m-1}(1))}{\sqrt{\delta(t)^{m-1}}}
\end{aligned} 
 \end{equation}
 and hence 

 \begin{equation}\label{final}
 \begin{aligned}
\frac{\Vol(\partial V(t))}{\Vol(S_t^w)}\leq &\frac{1}{w(t)^{m-1}\sqrt{\delta(t)^{m-1}}}\\&=\frac{1}{\left(1-4\epsilon(t)\right)^{(m-1)/2}}
\end{aligned} 
 \end{equation}
 Therefore, since for $t$ large enough $Vol(\partial D_t(o))\leq \sum_{i=1}^{\mathcal{E}(P)} Vol(\partial V_i(t))$ where $V_i$ denotes each end of $P$ then:
 \begin{equation}\label{final}
 \begin{aligned}
\frac{\Vol(\partial D_t(o))}{\Vol(S_t^w)}\leq \frac{\mathcal{E}(P)}{\left(1-4\epsilon(t)\right)^{(m-1)/2}}
\end{aligned} 
 \end{equation}
 

 \section{Proof of Theorem \ref{minimal}}\label{demos3}

To show assertion (1) we apply Theorem \ref{Isoperimetric} and inequality (\ref{lowerboundends}) in Theorem \ref{theorem2} to obtain, for $r$ sufficiently large, (we suppose that $\varphi^{-1}(o_w)\neq \emptyset$, and take $o \in \varphi^{-1}(o_w)$ in order to have that $\Vol(D_r(o)) \geq \Vol(B^w_r)$ for all $r >0$) :

\begin{equation}\label{final2}
\begin{aligned}
1\leq &\frac{\Vol(D_r(o))}{\Vol(B^w_r)} \leq \frac{\Vol(\partial D_r(o))}{\Vol(S_r^w)}\\ &\leq \frac{\mathcal{E}(P)}{\left(1-4\epsilon(r)\right)^{(m-1)/2}}
\end{aligned}
\end{equation}

Moreover, we know (again using Theorem \ref{Isoperimetric})  that the volume growth function is non-decreasing. 

Therefore, taking limits in (\ref{final2}) when $r$ goes to $\infty$, we obtain:

\begin{equation}\label{final3}
1\leq \lim_{r \to \infty}\frac{\Vol(D_r(o))}{\Vol(B^w_r)}=\Sup_{r>0} \frac{\Vol(D_r(o))}{\Vol(B^w_r)} \leq \mathcal{E}(P)
\end{equation}
Now, to prove assertion (2), we have, if $P$ has one end, that

\begin{equation}
1\leq \Sup_{r>0} \frac{\Vol(D_r(o))}{\Vol(B^w_r)}\leq 1
\end{equation}

Hence, as $f(r)=\frac{\Vol(D_r(o))}{\Vol(B^w_r)}$ is non-decreasing, then $f(r)=1 \,\,\forall r>0$, so we have equality in inequality (\ref{isopComp}) for all $r>0$, and $P$ is a minimal cone, (see \cite{MP} for details).



 \section{Proof of Theorems \ref{cor1} and \ref{cor3} and the Corollaries}\label{proofcorollaries}
 

\subsection{Proof of Theorem \ref{cor1}}
We are going to apply Theorem \ref{minimal}. To do that, we must to check hypotheses (1) and (2) in Theorem \ref{theorem2}. 

 We have, in this case, that the ambient manifold is the hyperbolic space $\Han$. Therefore all of its points are poles, so there exist at least $o \in \Han$ such that $\varphi^{-1}(o) \neq \emptyset$. As it is known, Hyperbolic space $\Han$ is a model space with $w(r)=w_b(r)=\frac{1}{\sqrt{-b}} \sinh \sqrt{-b} r$ so $w_b'(r) =\cosh\sqrt{-b} r \geq 1\,\, \forall r>0$. 
 
Therefore, hypothesis (2) in Theorem \ref{theorem2} is fulfilled in this context. Concerning hypothesis (1), it is straightforward that
\begin{equation}
\begin{aligned}
\Vert B^P_x \Vert &\leq \frac{\delta(r(x))}{e^{2\sqrt{-b}\, r(x)}}  \leq \frac{\epsilon(r)\sqrt{-b}}{\sinh\sqrt{-b} r \cosh\sqrt{-b} r}\\&= \frac{\epsilon(r)}{ \cosh^2\sqrt{-b} r}\sqrt{-b}\coth\sqrt{-b} r=\frac{\epsilon(r)}{(w_b'(r))^2}\eta_{w_b}(r)
\end{aligned}
\end{equation}
\noindent where $\epsilon(r)= \frac{\delta(r(x))}{4\sqrt{-b}}$ goes to $0$ when $r$ goes to $\infty$. 

Hence, also hypothesis (1) in Theorem \ref{theorem2} is fulfilled so, applying inequality (\ref{volgrowthends}) in Theorem \ref{minimal}, (because $P$ is minimal)
\begin{equation}\label{volgrowthends2}
1\leq \lim_{r\to \infty}\frac{Vol(D_r)}{Vol(B_r^{w_b})}\leq \mathcal{E}(P)
 \end{equation}

Finally, when $P$ has one end, then $ \lim_{r\to \infty}\frac{Vol(D_r)}{Vol(B_r^{w_b})}=1$. Since $P$ is minimal, by Theorem \ref{Isoperimetric}, $f(r)=\frac{Vol(D_r)}{Vol(B_r^{w_b})}$ is a monotone non-decreasing function, and, on the other hand, $f(r) \geq 1 \,\,\forall r>0$ because inequality (\ref{eqVolComp}). Hence $f(r)=1\,\,\forall r>0$, so $f'(r)=0\,\,\forall r>0$. This last equality implies the equality in inequality (\ref{isopComp}) for all $r>0$, (see \cite{MP} or \cite{MP1} for details), and we apply equality assertion in Theorem \ref{Isoperimetric} to conclude that $P$ is totally geodesic in $\Han$.
\bigskip


\subsection{Proof of Theorem \ref{cor3}}
In this case, we apply Theorem \ref{minimal}, being $M^n_w=\erre^n$, i.e., being $w(r)=w_0(r)=r$, ($b=0$). Hence, $w_0'(r)=1 >0 \,\,\forall r>0$ and $\eta_0(r) =\frac{1}{r}$ and hypotheses (1) and (2) in this theorem are trivially satisfied.
 
 When $P$ has only one end we conclude as before that the volume growth function is constant so we conclude equality in (\ref{isopComp})  for all radius $r>0$. Hence $P$ is totally geodesic in $\erre^n$ applying the corresponding equality assertion in Theorem \ref{Isoperimetric}.

\bigskip


\subsection{Proof of Corollary \ref{cor2}}

We are considering now a complete and proper immersion in $\Han$, as in Theorem \ref{cor1}, but $P$ is not necessarily minimal. In this setting hypotheses (1) and (2) in Theorem  \ref{theorem2} are fulfilled (as we have checked in the proof above, without using minimality). Hence taking limits in (\ref{lowerboundends}) when we consider an increasing sequence $\{t_i\}_{i=1}^\infty$ such that $t_i\to \infty$ when $i\to \infty$, we have:

$$\liminf_{i \to \infty}\frac{\Vol(\partial D_{t_i})}{\Vol(S_{t_i}^{b,m-1})}\leq \mathcal{E}(P)$$
\bigskip


\subsection{Proof of Corollary \ref{cor4}}
Hypotheses (1) and (2) in Theorem \ref{theorem2} are trivially satisfied and we argue as in the proof of Corollary \ref{cor2} to obtain the result.

\bigskip

\subsection{Proof of Corollary \ref{cor5}}
We apply Theorem \ref{theorem1}. Our ambient manifold is $\kan$, ($b \leq 0$), so hypothesis (1) about the bounds for the radial sectional curvature holds, and  as $w(r)=w_b(r)$ hence $w_b'(r)\geq 1 >0 \,\,\forall r>0$ and $\eta_{w_b}'(r) \leq 0 \,\,\forall r>0$. This means that hypothesis (3) is fulfilled. Hypothesis (2) in Theorem \ref{theorem1} holds because 
$$\Vert B^P_x \Vert \leq c\,h_b(\rho^P(x))$$
where $\rho^P(x)$ is the (intrinsic) distance to a fixed $x_o \in \varphi^{-1}(o)$ and  $c$ is a positive constant such that $c<1$.

\bigskip

\subsection{Proof of Corollary \ref{cor6}}
We apply again Theorem \ref{theorem1}, having into account that the ambient space is the Cartan-Hadamard manifold $\Han \times \erre^l$ and the model space used to compare is $\erre^m$, with $w(r)=w_0(r)=r$.

\bigskip


\end{document}